\documentclass[12pt]{amsart}
\usepackage{mathrsfs}
\usepackage{txfonts}
\usepackage{amssymb}

\makeatletter \@mparswitchfalse \makeatother
\def\eqref#1{(\ref{eq#1})}

\numberwithin{equation}{section}
\begin{document}
\title{Star  order automorphisms on the poset of type 1 operators }
\author{Xinhui Wang}\email{baichuan@snnu.edu.cn}\author{Guoxing Ji}\email{gxji@snnu.edu.cn}
  \address{School  of Mathematics and Information Science,
  Shaanxi Normal University,
  Xian , 710119, People's  Republic of  China}

\thanks{This research was
supported by the National Natural
   Science Foundation of China(No. 11771261) and the Fundamental Research Funds for the Central Universities (Grant No. GK201801011).}

    %\no{\footnotesize{\bf MSC(2000):\hspace{2mm}Primary    46L52, 47L75: Secondary 46K50,  46J15 }

   % \subjclass{     46L52\sep  46K50}
 \maketitle
\begin{abstract}Let $\mathcal{H}$ be a complex  infinite dimensional   Hilbert space and $\mathcal{B}(\mathcal{H})$ the algebra of all bounded linear operators on $\mathcal H$. The star partial order is defined by $A\overset{*}{\leq}B$ if and only if $A^*A=A^*B$ and $AA^*=AB^*$ for any  $A$ and $B$ in  $\mathcal B(\mathcal H)$. We give a type decomposition of operators with respect to star  order. For any $A\in\mathcal B(\mathcal H)$, there  are  unique  type 1  operator $A_1 $ which is $0$ or  the supremum of those rank 1 operators less than $A_1$    and type 2 operator $A_2$ which is not greater  than any  rank 1 operator  in  star order such that  $A_i\overset{*}{\leq}A$$(i=1,2)$  and $A=A_1+A_2$. Moreover, we determine all  automorphisms on the poset of type 1 operators.  As a consequence, we characterize continuous automorphisms on $\mathcal B(\mathcal H)$.
 \\
  {\bf keywords}{  star  order; type decomposition;   automorphism}\\
{\bf 2010 MSC}     06A06,    47B49
   \end{abstract}
\baselineskip18pt

\section{Introduction}

Partial order is a very important notion in many research areas.  In recent years, many   researchers make attention to star order on operator algebras.  Drazin in \cite{dr} introduced  the star order $\overset{*}{\leq}$ on matrix
algebra $M_{n}$ of all $n\times n$ complex matrices, that is,  for   $A, B\in M_n$, then we say that $A \overset{*}{\leq} B$  if  $A^{*}A=A^{*}B$ and $AA^{*}=BA^{*}$.  We note that  this definition can be extended to a $C^*$-algebra by the same way. In particular, it can be extended to the $C^*$-algebra
  $\mathcal B(\mathcal H)$ of all bounded linear operators on a complex infinite dimensional  Hilbert space $\mathcal H$. It is well-known that many very important structure results on this partial order are obtained(cf. \cite{an,de,do,ha,xu,zhang,zhangs}).
Among other things, characterizations of isomorphisms on certain partial orders are very interesting topics.
Guterman in \cite{gu}  characterized linear bijective maps on  $M_n$ preserving the star order  and Legi$\check{\mbox{s}}$a  in \cite{le} considered  automorphisms of $M_n$ with respect to the star order.
   Recently, several authors consider  star order preserving maps on  certain  subsets of $\mathcal B(\mathcal H)$ or a general von Neumann algebra  with respect to the star order when $\mathcal H$ is infinite dimensional. Dolinar and Guterman in \cite{dol} studied the automorphisms of the  algebra $\mathcal K(\mathcal H)$ of compact operators on a separable complex infinite dimensional  Hilbert space $\mathcal H$ and they  characterized the bijective,
additive, continuous maps on $\mathcal K(\mathcal H)$ which preserve the star order in both directions. An improvement of this result may found in \cite{xi}. On the other hand,
 characterizations of certain continuous  bijections on the normal elements of a von Neumann algebra  preserving the star order in both directions  are considered(cf.\cite{bo1,bo2}). However, star order automorphisms on $\mathcal B(\mathcal H)$ are still unknown. For example, what is an infinite dimensional version of Legi$\check{\mbox{s}}$a's result  in \cite{le}?
   In this paper, we consider this problem. We consider a type decomposition of operators with respect to star order.
    We then determine star order automorphisms on the poset of type 1 operators. As a consequence, we determine all  continuous star order automorphisms  on $\mathcal B(\mathcal H)$.

Let ${\mathcal H}$  be    a complex  infinite dimensional
 Hilbert
 space and let   ${\mathcal B({\mathcal H})}$, $\mathcal K(\mathcal H)$ and $\mathcal F(\mathcal H)$
be the  algebra   of all bounded linear operators, the compact operators and the finite rank operators  on  ${\mathcal H}$, respectively.
For every pair of vectors $x,
y\in\mathcal H,$  $\langle x,y \rangle  $ denotes the inner
product of $x$ and $y$, and
$x\otimes y$ stands for    the   rank-1  linear operator on $\mathcal H$
defined by $(x\otimes y)z=\langle z, y\rangle x$ for any
$z\in\mathcal H$.  If $x$ is a unit vector, then $x\otimes x$ is a rank-1 projection. $\sigma(A)$  and $\sigma_p(A)$  are   the spectrum and point spectrum
of $A$   for any $A\in\mathcal B(\mathcal H)$.  For
  a subset  $S$ of $\mathcal H$,  $[S]$ denotes the closed subspace of
  $\mathcal H$ spanned by $S$ and $P_M$ denotes the orthogonal projection on $M$ for a closed subspace  $M$ of $\mathcal H$.
 We denote by
$R(T)$ and  $\ker T$  the range and  the kernel  of a linear map $T$  between two linear spaces.   Throughout this paper, we will generally  denote by $I $ the
identity operator on  a Hilbert space.
\section{A type  decomposition  of  operators with respect to star order}

Let $A, B\in\mathcal B(\mathcal H)$. We recall that  $A \overset{*}{\leq} B$  if  $A^{*}A=A^{*}B$ and $AA^{*}=BA^{*}$. This partial order is called star order.
 Let  $H_1 $, $H_2$, $K_1$    and $K_2 $ be closed subspaces of $\mathcal H$ such that
 $$\mathcal H=H_{1}\oplus H_{2}=K_{1}\oplus K_{2}.  \eqno(2.1)$$
 Let $A,\ B\in\mathcal B(\mathcal H)$ with   the following matrix forms
$$A=\left(
\begin{array}{ccc}
A_{11}& A_{12} \\A_{21} & A_{22} \\
\end{array}
\right)  \mbox{ and } B=\left(
\begin{array}{ccc}
B_{11}& B_{12} \\B_{21}& B_{22} \\
\end{array}
\right)\eqno(2.2)$$ with respect to the orthogonal  decompositions $(2.1)$, where $A_{ij},B_{ij}\in\mathcal B(H_j,K_i)$ $(i,j=1,2)$.
 It is known that $A\overset{*}{\leq}B$ if and only if there exists a decomposition $(2.1)$ of $\mathcal H$ such that $A_{11}=B_{11}$ and   $A_{12}$,
 $A_{21}$, $A_{22}$ as well as $B_{12}$, $B_{21}$ are 0(cf.\cite{do}). In fact, we may choose
  $H_1=\overline{R(A^*)}$, $H_2=\ker A$, $K_1=\overline{R(A)}$ and $K_2=\ker A^*$ respectively.

  We say that $A$ and $B$ are orthogonal and denoted by $A\bot B$ if $A^*B=AB^*=0$.    Note that $A\bot B$ if and only if there is a decomposition  $(2.1)$  such that  all entries of $A$ and $B$  in  matrix forms $(2.2)$ are $0$ except $A_{11}$ and $B_{22}$.  The following result are elementary.
 \vskip12pt

 {\bf Proposition 2.1} Let $A,\  B\in \mathcal B(\mathcal H)$. Then $A\overset{*}{\leq }A+B$(resp. $B\overset{*}{\leq}A+B$) if and only if $A\bot B$.
 \vskip12pt

 Let $\mathcal{A}\subseteq\mathcal{B(H)}$. If there exists an operator $B\in\mathcal{B(H)}$ such that $A\overset{*}{\leq} B$, for all $A\in\mathcal{A}$, then we say $B$ is an upper bound of $\mathcal A$. If $D$ is an upper bound of $\mathcal A$
such that for  any  upper  bound $B$ of  $\mathcal A$, we have
 $D\overset{*}{\leq}B$,   then we say that $D$ is   the supremum of $\mathcal A$ and  denoted by $\sup\mathcal A$.  It is known that a subset  with an upper bound has the supreum(cf.\cite{xu,zhang}).

 \vskip12pt
{\bf Definition 2.2} Let $A\in \mathcal B(\mathcal H)$. If $A=0$ or  $A= \sup\{x\otimes y:\ x\otimes y \overset{*}{\leq} A\}$, then we say that $A$ is of type 1 with respect to star order. If there are no rank 1 operators $x\otimes y$ such that $x\otimes y\overset {*}{\leq}A$, then we say that $A$ is of type 2 with respect to star order.
\vskip12pt

Put $\mathcal B(\mathcal H)_I=\{A\in\mathcal B(\mathcal H): A \mbox{ is  of type 1 }\}$ and $\mathcal B(\mathcal H)_{II}=\{A\in\mathcal B(\mathcal H): A \mbox{ is of  type 2}\}$. It is elementary both $\mathcal B(\mathcal H)_{I}$ and $\mathcal B(\mathcal H)_{II}$ are posets with respect to star order. By use of polar decomposition as well as spectral decomposition, we easily have

\vskip12pt
{\bf Proposition 2.3} Let $A,\ B\in\mathcal B(\mathcal H)$.

$(1)$ If $A\overset{*}{\leq}B$ and $B$ is of  type $i$ $(i=1,2)$, then so is $A$.

$(2)$
$\mathcal B(\mathcal H)_I$ is dense in $\mathcal B(\mathcal H)$ in norm topology.
\vskip12pt

However, $\mathcal B(\mathcal H)_{II}$ is not dense in norm topology. In fact,
For any  nonzero $K\in\mathcal K(\mathcal H)$, $K \not\in\overline{\mathcal B(\mathcal H)_{II}}$. Otherwise, if $\|A_n-K\|\to 0$$(n\to\infty)$ for a sequence $\{A_n:n=1,2,  \cdots \}\subseteq {\mathcal B(\mathcal H)_{II}}$, then $\|A_n\|_e=\|A_n-K\|_e\to 0$$(n\to\infty)$.  However, it is easy to show that $\|A\|=\|A\|_e$ for any $A\in\mathcal B(\mathcal H)_{II}$, where $\|A\|_e$ is the essential norm of $A$. Then
 $\|A_n\|=\|A_n\|_e\to 0$$(n\to\infty)$. Note that $\|A_n-K\|>\|K\|-\|A_n\|_e$ for large $n$. This is a  a contradiction. Therefore $\mathcal B(\mathcal H)_{II}$ is not dense in norm topology.

\vskip12pt
{\bf Theorem  2.4} Let $A\in \mathcal B(\mathcal H)$.

$(1)$ There exists unique pair of operators $A_1\in\mathcal B(\mathcal H)_I$ and $A_2\in \mathcal B (\mathcal H)_{II}$ such that $A_i\overset{*}{\leq}A$ for  $i=1,2$ and $A=A_1+ A_2$.

$(2)$ $A\in \mathcal B(\mathcal H)_I$ if and only if there exist a  unique family of   mutually orthogonal nonzero partial isometries $\{U_j:j\in J\}$ and a family  of  positive numbers $\{a_j: j\in J\}$ with $a_i\not=a_j$ for any $i,j\in J,\ i\not=j$,   such that
 $$ A=\sum\limits_{j\in J}a_jU_j.\eqno(2.3)$$
\vskip12pt

\begin{proof} $(1)$ If there are  no rank 1 operators $x\otimes y$ such that $x\otimes y\overset {*}{\leq}A$, then  $A_1=0$ and $A_2=A$.

Assume that $\mathcal A_1=\{x\otimes y:\ x\otimes y \overset{*}{\leq} A\}\not=\emptyset$. Then $\mathcal A_1 $ is a subset with upper bound $A$. Put $A_1=\sup\mathcal A_1$. Then $A_1$ is of type 1 and $A_1\overset{*}{\leq}A$. It is well known that $A_2=A-A_1\overset{*}{\leq}A$  with  $A_1 \bot A_2$. If there is a rank 1 operator $x\otimes y\overset{*}{\leq}A_2$, then $x\otimes y\overset{*}{\leq}A_2\overset{*}{\leq }A$. Thus $x\otimes y\overset{*}{\leq}A_1$. This is a contradiction. Therefore $A_2$ is of type 2 and $A=A_1+A_2$.  Assume that  $A=B_1+ B_2$ with $B_1\in \mathcal B(\mathcal H)_1$ and $B_2\in\mathcal B(\mathcal H)_{II}$ such that $B_i\overset{*}{\leq}A$$(i=1,2)$.  Then for any rank 1 operator $x\otimes y\overset{*}{\leq}B_1\overset{*}{\leq}A$, we have that $x\otimes y\overset{*}{\leq}A_1$. It follows that $B_1\overset{*}{\leq}A_1$.
Note that $A_1-B_1\overset{*}{\leq}A_1$ and  $A_1=B_1+ (A_1-B_1)$. Thus   $A_2 \bot (A_1-B_1)$ and $B_2=A_2 +(A_1-B_1)$. Note  that $A_1-B_1$ is also  of type 1 by Proposition 2.3. There is a rank 1 operator $x\otimes y\overset{*}{\leq}A_1-B_1\overset{*}{\leq}B_2$, a contradiction. Therefore $B_1=A_1$ and $B_2=A_2$.

$(2)$ Let $A\in\mathcal B(\mathcal H)_I$ and $\mathcal A_1$ as above. We take a maximal mutually orthogonal family $\{x_i\otimes y_i:   i\in \Lambda \}$ in $\mathcal A_1$ by Zorn's lemma. It is trivial that
$A=\sum_{i\in\Lambda } x_i\otimes y_i$.  Now we define an equivalence $\sim$ in this family by $x_i\otimes y_i\sim x_l\otimes y_l$ if and only if $\|x_i\|\|y_i\|=\|x_l\|\|y_l\|$ and put $\Lambda_j=\{l: x_l\otimes y_l\sim x_j\otimes y_j\}$ is the equivalent class for any $j$   in  a set $J$.
       We now define
$a_j=\|x_i\|\|y_i\|$ for some $i\in \Lambda_j$ and $U_j=\sum_{i\in {\Lambda_j}}\frac{1}{a_j}x_i\otimes y_i.$ Then $U_j$ is a partial isometry  with initial space $\vee\{y_i:i\in \Lambda_j\}$ as well as final space $\vee\{x_i:i\in \Lambda_j\}$ for every $j$ and     $A=\sum\limits_{j\in J} a_j U_j$.  Put $E_j=U_j^*U_j$$(1\leq j\leq k)$ and $U=\sum\limits_{j\in J}U_j$. It is known that $|A|=\sum_{j\in J}a_jE_j$ and $\sigma_p(|A|)-\{0\}=\{a_j: j\in J\}$.     Note that  $A=U|A|$ is the polar decomposition of $A$.

 On  the other hand, put   $A=\sum\limits_{i\in K}b_iV_i$ for some mutually different positive numbers $\{b_i: i\in K\}$ and mutually orthogonal partial isometries
 $\{V_i: i\in K\}$. Then  $|A|=\sum\limits_{i\in K}b_iV_i^*V_i$. Thus $\sigma_p(|A|)-\{0\}=\{b_i:i\in K\}$. This means that $\{a_j:j\in J\}=\{b_i:i\in K\}$ and therefore the cardinalities of $J$ and $K$ are the same, $a_j=b_j$ and $E_j=V_j^*V_j$ by rearranging  $\{b_i:i\in K\}$.  It follows that $V_j=U_j$  for any $j\in J$.

  The converse is clear.
\end{proof}

 We call  $A=A_1+A_2$ as in Theorem 2.4  the type decomposition of $A$ with respect to  star order. As in \cite[Theorem 3.1]{le},  we call $(2.3)$ the Penrose  decomposition of $A$.  The following theorem is similar to \cite[Theorem 3.3]{le}.

\vskip12pt
{\bf Theorem 2.5}
If $A, \ B\in\mathcal B(\mathcal H)_I$ have the Penrose decompositions
$$ A=\sum\limits_{j\in J}a_jU_j, \ B=\sum\limits_{i\in K}b_iV_i, $$ then
$A\overset{*}{\leq }B$ if and only if there is an injection $\Lambda:\ J\to K$ such that
$a_j=b_{\Lambda(j)}$ and $U_j\overset{*}{\leq} V_{\Lambda(j)}$ for all $j\in K$.

\section{Automorphisms of $\mathcal B(\mathcal H)_I$}
Let $\varphi$ be a map on $\mathcal B(\mathcal H)$. We say that $\varphi$ is a star order automorphism if  $\varphi$ is bijective  such that $\varphi(A)\overset{*}{\leq}\varphi(B)\Leftrightarrow A\overset{*}{\leq}B$, $\forall A, \ B\in\mathcal B(\mathcal H)$.  The following proposition is elementary.

\vskip12pt
{\bf Proposition 3.1} Let  $\varphi$ be  a star order automorphism on $\mathcal B(\mathcal H)$. Then for all  $A\in \mathcal B(\mathcal H)$, $\mbox{rank}\varphi(A)=\mbox{rank} A=\mbox{rank}\varphi^{-1}(A)$. In particular, $\varphi(0)=0$ and $\varphi(A)$ is of rank 1 if and only is $A$ is.
 \vskip12pt

{\bf Theorem 3.2} Let $\varphi$ be a bijection on $\mathcal B(\mathcal H)$. then $\varphi$ is a  star order  isomorphism  if and only if there exist star order automorphisms $\varphi_1$ and $\varphi_2$ on  $\mathcal B(\mathcal H)_I$ and $\mathcal B(\mathcal H)_{II}$  such that $\varphi(A)=\varphi_1(A_1)+\varphi_2(A_2)$  for all $A=A_1+A_2$,  the type   decomposition of $A$,  where $A_1\in \mathcal B(\mathcal H)_I$ and $A_2\in \mathcal B(\mathcal H)_{II}$.
\vskip12pt

\begin{proof}
Let $\varphi $ be a star order automorphism on $\mathcal B(\mathcal H)$. Then it is elementary that $\varphi(A)$ is of type $i$($i=1,2$) if and only if $A$ is by Proposition 3.1. That is,  $\varphi_1=\varphi|_{\mathcal B(\mathcal H)_I}$ and
$\varphi_2=\varphi|_{\mathcal B(\mathcal H)_{II}}$ are automorphisms. Thus for any $A=A_1+A_2$, the type  decomposition of $A\in\mathcal B(\mathcal H), $
$\varphi(A)=\varphi(A_1)+\varphi(A_2)=\varphi_1(A_1)+\varphi_2(A_2)$.

The converse is clear.
\end{proof}
  We next give an example of  star partial automorphism on $\mathcal B(\mathcal H)$. For any $A\in\mathcal B(\mathcal H)$, we denote by $A=W_A|A|$ the polar decomposition of $A$.

\vskip12pt
{\bf Example 3.3}  Let $f$  be a  bijection on $(0,\infty)$ such that $f$ is bounded on any bounded subset of $(0,\infty)$ and $g$ a continuous bijection on $[0,\infty)$ with $g(0)=0$.  Let $\varphi_1(0)=0$. For any nonzero $A=\sum\limits_{j\in J}a_jU_j\in \mathcal B(\mathcal H)_I$, we define $\varphi_1(A)=\sum\limits_{j\in J}f(a_j)U_j$. For any $B\in\mathcal B(\mathcal H)_{II}$,
we define $\varphi_2(B)=W_Bg(|B|)$. If $\varphi(A)=\varphi_1(A_1)+\varphi_2(A_2)$, for all $A\in\mathcal B(\mathcal H)$, where $A=A_1+A_2$ is the type decomposition  with respect to star order, then $\varphi$ is a star order automorphism on $\mathcal B(\mathcal H)$.
\vskip12pt

  By Theorem 3.2, it is sufficient  to  consider separately star order automorphisms on $\mathcal B(\mathcal H)_I$ and $\mathcal B(\mathcal H)_{II}$ respectively. We next assume that $\varphi$ is a star order automorphism on $\mathcal B(\mathcal H)_I$.   Put
  $$PI(\mathcal H)=\{V\in\mathcal B(\mathcal H):  V\mbox{  is a nonzero partial isometry} \}.$$

\vskip12pt
{\bf Lemma 3.4} Let $\varphi$ be  a star order  isomorphism on $\mathcal B(\mathcal H)_I$. Then for any  constant $a \in (0,\infty)$ and
 $U\in PI(\mathcal H)$, there are   a  constant $f(a ,U)\in  (0,\infty)$, a  partial isometry
 $W(a, U)\in PI(\mathcal H)$ such that  $\varphi(a  U)=f(a ,U)W(a, U)$. If $U$ is unitary, so is $W(a, U)$.
\vskip12pt

 \begin{proof} We  firstly  assume that $a =1$. It is clear that $V$ is a partial isometry for any $V\overset{*}{\leq}U$. Put $B=\varphi(U)$ and $B=W_B|B|$ is the polar decomposition of $B$. We claim that there is at most one nonzero number in $\sigma(|B|)$. Otherwise, let $ \alpha,\beta\in\sigma(|B|)$ be two nonzero positive numbers and   $|B|=\int \lambda dE_{\lambda}$ is the spectral decomposition of $|B|$. Take any two closed intervals $[\alpha_1,\alpha_2]$ and $[\beta_1,\beta_2]$ such that $\alpha_1>0$, $\beta_1>0$,  $[\alpha_1,\alpha_2]\cap [\beta_1,\beta_2]=\emptyset$, $\alpha\in [\alpha_1,\alpha_2]$ and $\beta\in [\beta_1,\beta_2]$.
  Let $\mathcal H=H_1\oplus H_2\oplus H_3$, where $H_1=E[\alpha_1,\alpha_2]\mathcal H$, $H_2= E[\beta_1,\beta_2]\mathcal H$ and
   $H_3=\mathcal H\ominus (E[\alpha_1,\alpha_2]\mathcal H\oplus E[\beta_1,\beta_2]\mathcal H)$. Note that $P_i=|B||_{H_i}$ is invertible and  $W_i=W|_{H_i}$ is unitary for $i=1,2$. Take two rank 1 operators $x_i\otimes y_i$ such that $x_i\otimes y_i\overset{*}{\leq}W_iP_i$ for $i=1,2$.  Then $B_2=x_1\otimes y_1 +x_2\otimes y_2 \overset{*}{\leq} W_1P_1 +W_2P_2 \overset{*}{\leq}B$.  It is elementary that  there are only two rank 1 operators $x_i\otimes y_i \overset{*}{\leq}B_2$ for $i=1,2$.
    Put $V_2=\varphi^{-1}(B_2)$. Note that  $V_2\overset{*}{\leq}U$ is a rank 2 partial isometry by Proposition 3.1.  Thus for any unit vector
    $x$ in the initial space of $V_2$, we have $V_2x\otimes x\overset{*}{\leq}V_2$ and $\varphi(V_2x\otimes x)\overset{*}{\leq} B_2$.
    This is a contradiction. Therefore  there is at most one nonzero number $b$  in $\sigma(|B|)$ and $|B|=bE$ for a projection $E$. This implies that
   $B=bW_B$.

    In general, we define $\psi(T)=\varphi(aT)$, $\forall T\in\mathcal B(\mathcal H)_I$. Then $\psi$ is also a star order automorphism on $\mathcal B(\mathcal H)_I$ with $\psi(U)=\varphi(aU)$. By  the above proof, $\varphi(aU)=\psi(U)=bW $  for  a nonzero constant $b$ and  a partial isometry $W$. We  may define $f(a,U)=b$ and $W(a,U)=W$.

  We now assume that $U$ is unitary and $\varphi(U)=bW$. Let $\psi(T)=b^{-1}\varphi(UT)$, $\forall T\in\mathcal B(\mathcal H)$. Then $\psi$ is also a star order isomorphism such that $\psi(I)=b^{-1}\varphi(U)=W$. Without loss of generality,  we show this for $U=I$ and $\varphi(I)=W$. We firstly claim that  $W$ is an isometry or a co-isometry. Otherwise, there is a partial isometry $W_0$ such that $W\overset{*}{\leq}W_0$ and $W\not=W_0$. Then $I\overset{*}{\leq}\varphi^{-1}(W_0)$ and $I\not=\varphi^{-1}(W_0)$, a contradiction.  We next   assume that $W$ is a non-unitary  isometry.  That is, $\ker W^*\not=\{0\}$.

  For any  unit vector $x$, put $E_x=I-x\otimes x $ and $W_x=\varphi(E_x)$.  Then $W_x\overset{*}{\leq}W$ and there is a rank 1 partial isometry $z_x\otimes y_x$  with $W_x\bot z_x\otimes y_x$ such that $W=W_x + z_x\otimes y_x$.  Take any unit vector $z\in\ker W^*$. Note that $W_x\bot \mu z\otimes y_x$ for all nonzero $\mu\in\mathbb C$. Put $W_{x,\mu z}=W_{x}  +\mu z\otimes y_x$ for any nonzero   $\mu\in\mathbb C$. It is clear that $W_{x,\mu z}\not=W_x$ and $W_x\overset{*}{\leq}W_{x,\mu z}$. Then $E_x\not=\varphi^{-1}(W_{x,\mu z})$ and  $E_x\overset{*}{\leq}\varphi^{-1}(W_{x,\mu z})$. In fact, we note that $\{T\in\mathcal B(\mathcal H):E_x\overset{*}{\leq}T\}=\{E_x + \gamma x\otimes x: \gamma\in\mathbb C\}$. Thus there is a nonzero $a(x,\mu, z)\in\mathbb C$ such that $\varphi^{-1}(W_{x,\mu z})=E_x + a(x,\mu,z)x\otimes x$
for any nonzero $\mu\in\mathbb C$. Since $W_{x,\mu z}$ is not a partial isometry for any nonzero  $|\mu|\not=1$,
$|g(x,\mu,z)|\not=1$.
It is easy to know that $a(x,\mu,z)x\otimes x$ is the unique rank 1 operator  such that $a(x,\mu,z)x \otimes x\overset{*}{\leq}E_x + a(x,\mu,z)x\otimes x$ for $|a(x,\mu,z)|\not=1$ by Theorem 2.5. So is $\mu z\otimes y_z$ such that
 $\mu z\otimes y_z\overset{*}{\leq} W_x + \mu z\otimes y_x$ for  nonzero $|\mu|\not=1$.  This means that $\varphi(a(x,\mu,z)x\otimes x)=\mu z\otimes y_x$.

 Take two orthogonal unit   vectors $x_i(i=1,2)$ and put $A=a(x_1,\mu,z)x_1\otimes x_1+a(x_2,\mu, z)x_2\otimes x_2$
 for  nonzero $|\mu|\not=1$ and $B=\varphi(A)$. Then $A$ is of rank 2 and  $a(x_i,\mu,z)x_i\otimes x_i\overset{*}{\leq}A$
 for $1=1,2$. Thus $\mu z\otimes y_{x_i}\overset{*}{\leq} B$ for $i=1,2$. Note that $y_{x_1}\not=y_{x_2}$. This is a contradiction. Therefore $\ker W^*=\{0\}$ and  $W$ is unitary.
 \end{proof}

  Let $\varphi(I)=f(1,I)W(1,I)$. We again define $\phi(T)=f(1,I)^{-1}W(1,I)^*\varphi(T)$, $\forall T\in\mathcal B(\mathcal H)_I$. Then $\phi$ is a star order automorphism such that $\phi(I)=I$.    We next assume that $\varphi(I)=I$.

\vskip12pt
 {\bf Lemma 3.5}  Let $\varphi$ be  a star order  isomorphism on $\mathcal B(\mathcal H)_I$ with  $\varphi(I)=I$. Then $f(a, U)=f(a,I) $ for all $a \in (0,\infty)$, $U\in PI(\mathcal H)$. Moreover,  $\varphi( P I(\mathcal H))= P I(\mathcal H)$.
\vskip12pt

\begin{proof} Assume $a\in (0,\infty )$. Note that $f(a, V)=f(a, U)$ for any $V\overset{*}{\leq}U$.
 Let $U\in PI(\mathcal H)$.    Take a rank 1 partial isometry $x\otimes y \overset{*}{\leq}U$.  Then  $f(a, x\otimes y)=f(a, U)$. Take a unit vector $z\in\{x,y\}^{\bot}$ and put $A=x\otimes y +z\otimes z$. Then $A$ is a partial isometry. Note that $x\otimes y\overset{*}{\leq}A$ and
  $z\otimes z\overset{*}{\leq}A$. Then  $f(a, U)=f(a, x\otimes y)=f(a, A)=f(a z\otimes z)=f(a,I)$.
  Since $f(1,I)=1$, $\varphi(PI(\mathcal H))=PI(\mathcal H)$.\end{proof}

    We now define $f(a)=f(a, I)$, $\forall a\in(0,\infty)$. it is trivial that $f$ is bijective on $(0,\infty)$.

\vskip12pt
{\bf Lemma 3.6} Let $\varphi$ be  a star order  isomorphism on $\mathcal B(\mathcal H)_I$ with  $\varphi(I)=I$. Then there is a unitary or an anti-unitary operator  $U$  such that $\varphi(E)=UEU^*$ for all projection $E$.
\vskip12pt

\begin{proof}
Note that   $\varphi(E)$ is a projection if and only if  $E$ is since  $\varphi(I)=I$.  Thus there exists a bounded invertible linear or conjugate linear operator $S$ on $\mathcal H$ such that $\varphi(E)=P_{SE}$, where $P_{SE}$ is the projection on the subspace $S(E\mathcal H)$ by \cite[Theorem 1]{fi}.

  For any unit vector $x$, as  in the proof of Lemma 3.4, Put $E_x=I-x\otimes x=P_{\{x\}^{\bot}}$. Then $\varphi(E_x)=P_{S\{x\}^{\bot}}$. Note that $S\{x\}^{\bot}$ is  of co-dimension 1.  For any $\lambda \in\mathbb C-\{0\}$ such that $|\lambda|\not=1$, put $A(x,\lambda)=E_x +\lambda x\otimes x$ and  $\varphi(A(x,\lambda))=B(x,\lambda)$. Since $E_x\overset{*}{\leq}A(\lambda,x)$, $P_{S\{x\}^{\bot}}\overset{*}{\leq}B(\lambda,x)$. Take a unit vector $\xi_x\in [S\{x\}^{\bot}]^{\bot}$. Then $B(\lambda,x)=P_{S\{x\}^{\bot}}\oplus g(\lambda, x)\xi_x\otimes \xi_x$ for a constant $g(\lambda, x)\in\mathbb C-\{0\}$  with $|g(\lambda, x)|\not=1$  by Lemma 3.4 since $A(\lambda,x)$ is not unitary. Note that if $h\otimes k\overset{*}{\leq}A_x$ is of rank 1, then either $h\otimes k\overset{*}{\leq}E_x$ or $h\otimes k=\lambda x\otimes x$  by Theorem 2.5. It now follows that $\varphi(\lambda x\otimes x)=g(\lambda, x)\xi_x\otimes \xi_x=f(|\lambda|)W(|\lambda|,\frac{\lambda}{|\lambda|}x\otimes x)$ by Theorem 2.5, Lemmas  3.4 and 3.5. Thus $|g(\lambda,x)|=f(|\lambda|)$  and  $\frac{g(\lambda,x)}{f(|\lambda|)}\xi_x\otimes \xi_x=W(|\lambda|,\frac{\lambda}{|\lambda|}x\otimes x)$.

   Now take any two orthogonal unit vectors $e_i$ and distinct nonzero numbers $\lambda_i$  with $|\lambda_i|\not=1$ $(i=1,2)$.
   Note that $\lambda_i e_i\otimes e_i\overset{*}{\leq} \lambda_1e_1\otimes e_1+\lambda_2 e_2\otimes e_2)$ for $i=1,2$ are the only two rank 1 operators with this property. Then $\varphi(\lambda_1e_1\otimes e_1+\lambda_2 e_2\otimes e_2)$ is of rank 2 such that
  $g(\lambda_i, e_i)\xi_{e_i}\otimes \xi_{e_i}\overset{*}{\leq} \varphi(\lambda_1e_1\otimes e_1+\lambda_2 e_2\otimes e_2)$  for $i=1,2$ are the only two rank 1 operators with this property.  This implies that $|g(\lambda,e_1)|=f(|\lambda_1|)\not=f(|\lambda_2|)=|g(\lambda_2,e_2)|$ and $\xi_{e_1}\bot \xi_{e_2}$.

     We extend $\{e_1,e_2\}$ to an orthogonal basis $\{e_{\gamma}:\gamma\in\Gamma\}$  of $\mathcal H$   and fixed a nonzero $\lambda$ with $|\lambda|\not=1$.
  Now $\lambda I=\sup \{\lambda e_{\gamma}\otimes e_{\gamma}:\gamma\in\Gamma\}$. Then
   $$f(|\lambda|) W(|\lambda|,\frac{\lambda}{|\lambda|}I)=\varphi(\lambda I)=\sup\{g(\lambda,e_{\gamma})\xi_{e_{\gamma}}\otimes\xi_{e_{\gamma}}:\gamma\in \Gamma\}$$ is  multiple of a unitary operator. Note that $\{\xi_{e_{\gamma}}:\gamma\in\Gamma\}$ is an  orthogonal unit  family in $\mathcal H$. This implies that
     $\{\xi_{e_{\gamma}}:\gamma\in\Gamma\}$ is an  orthogonal basis of $\mathcal H$.
      Since $\xi_{e_{\gamma}}\in [S(\{e_{\gamma}\}^{\bot})]^{\bot}$,
      $\{\xi_{e_{\nu}}: \nu\not=\gamma\} $ is an orthogonal basis of
   $ S(\{e_{\gamma}\}^{\bot})$.
    This implies that $Se_{\gamma}\in\cap_{\nu\not=\gamma}S(\{e_{\nu}\}^{\bot})=\vee\{\xi_{e_{\gamma}}\}$.
 In particular,  $Se_1\bot Se_2$. It follows that there exist a constant $c\in\mathbb C$ and   a unitary or an anti-unitary operator $U$ on $\mathcal H$ such that $S=cU$. Note that  $\varphi(E)=P_{SE(\mathcal H)}=P_{UE(\mathcal H)}=UEU^*$ for all projection $E$.
 \end{proof}

 Now we give the main result of this section.

 \vskip12pt
 {\bf Theorem 3.7} Let $\varphi$ be a star order automorphism on $\mathcal B(\mathcal H)_I$. Then $\varphi(0)=0$ and  there exist a  complex function $h$  on  $(0,\infty) $  bounded on any bounded subset so that $|h|$ is  bijective on $(0,\infty) $, a nonzero constant $\alpha\in (0,\infty)$ and  two unitary operators or two anti-unitary operators $S$ and $T$ on $\mathcal H$, such that one of the following assertions holds.

    $(1)$ $\varphi(A)=\alpha S(\sum\limits_{j\in J}h(a_j)U_j)T$ for all $A=\sum\limits_{j\in J}a_jU_j\in\mathcal B(\mathcal H)_I$;

    $(2)$ $\varphi(A)=\alpha S(\sum\limits_{j\in J}h(a_j)U_j^*)T$ for all $A=\sum\limits_{j\in J}a_jU_j\in\mathcal B(\mathcal H)_I$.
   \vskip12pt

   \begin{proof} Let $\varphi(I)=\alpha W$ for a positive constant $\alpha$ and a unitary operator $W$ by Lemma 3.4. Define
   $\phi(A)=\alpha^{-1}W^*\varphi(A)$ for all $A \in \mathcal B(\mathcal H)_I$. Then $\phi$ is a star order automorphism such
   that  $\phi(I)=I$. By Lemma 3.5, there exist
 a unitary or an anti-unitary operator $U$ on $\mathcal H$ such that $\phi(E)=UEU^*$ for all projection $E$.
 Let  $\psi(A)=U^*\phi(A)U$,$\forall A\in \mathcal B(\mathcal H)_I$. Then $\psi$ is a star order automorphism  such that $\psi(E)=E$ for all projection $E$.

 {\bf Claim 1 }  There is a function $g$  on  $\mathbb C-\{0\}$   such that $|g(\lambda)|=f(|\lambda|)$ and  $\psi(\lambda E)=g(\lambda) E$ for all $\lambda\in\mathbb C-\{0\}$ and projections $E\in\mathcal B(\mathcal H)_I$.

  We firstly assume that $|\lambda |\not=1$. By Lemma 3.6, we have $\psi(\lambda x\times x)=g(\lambda, x)x\otimes x$ for any rank 1 projection with $|g(\lambda,x)|=f(|\lambda|)\not=1$. Let  $E=\sum\limits_{j=1}^kx_j\otimes x_j$ be a rank k projection.
   Note that $\lambda E=\sup\{\lambda x_j\otimes x_j: 1\leq j\leq k\}$. Then
    $\psi(\lambda E)=\sup\{\psi(\lambda x_j\otimes x_j): 1\leq j\leq k\}=\sum\limits_{j=1}^kg(\lambda,x_j)x_j\otimes x_j$.
     For any unit vector $x\in\vee\{x_j:1\leq j\leq k\}$,
     $g(\lambda, x)x\otimes x \overset{*}{\leq}\sum\limits_{j=1}^k g(\lambda,x_j)x_j\otimes x_j$
     since $\lambda x\otimes x\overset{*}{\leq}\lambda E$. It easily follows that $g(\lambda,x\otimes x)=g(\lambda, x_j)$$(1\leq j\leq k)$.  This means that $g(\lambda, x)=g(\lambda, y)$ for all unit vectors $x,y \in\mathcal H$, that is, $g(\lambda,x)=g(\lambda)$ independent of $x$ and $\psi(\lambda E)=g(\lambda)E$ for any projection $E$.

      Let $|\lambda|=1$. Note that $\lambda x\otimes x $ is a partial isometry. So is $\psi(\lambda x\otimes x)$.
     Moreover, $ 2E_x +\lambda x\otimes x=\sup \{2E_x, \lambda x\otimes x \}$ for any unit vector $x$. Then
       $\psi(2E_x +\lambda x\otimes x)=\sup\{g(2)E_x, \psi(\lambda x\otimes x)\}$. It follows that $\psi(2E_x +\lambda x\otimes x)=g(2)E_x+g(\lambda,x)x\otimes x$ for a constant $g(\lambda,x)$. Similarly we have $|g(\lambda,x)|=1$ and
        $\psi(\lambda x\otimes x)=g(\lambda,x) x\otimes x$.  By use of the same method as above, we have that $g(\lambda,x)=g(\lambda) $ independent of $x$ and $\psi(\lambda E)=g(\lambda )E$ for any projection $E$.

{\bf Claim 2} For any $A=\sum\limits_{j\in J} a_j U_j\in\mathcal B(\mathcal H)_I$, $\psi(a_jU_j)=f(a_j)W(a_j,U_j)$  for some partial isometries $W(a_j,U_j)$$( j\in J)$ with
$W(a_i,U_i)\bot W(a_j,U_j)$ for all $i\not=j$ and $\psi(A)=\sum\limits_{j\in J} f(a_k)W(a_j, U_j)$.

 By Lemmas 3.4 and 3.5, we have $\psi(a_jU_j)=f(a_j)W(a_j,U_j)$ for some partial isometries $W(a_j,U_j)$($1\leq j\leq k)$.
  We note that for any  orthogonal rank 1 partial isometries  $x_i\otimes y_i$$(i=1,2)$,  $a_ix_i\otimes y_i$ for $i=1,2$  are
only two rank 1 operators such that $a_ix_i\otimes y_i\overset{*}{\leq}a_1x_1\otimes y_1+a_2x_i\otimes y_2$. This implies that
so are $\psi(a_ix_i\otimes y_i)$ for $i=1,2$ with
 $\psi(a_ix_i\otimes y_i)\overset{*}{\leq}\psi(a_1x_1\otimes y_1+a_2x_i\otimes y_2)$. It follows that
$\psi(a_1x_1\otimes y_1)\bot \psi(a_1x_2\otimes y_2)$  and $\psi(a_1x_1\otimes y_1+a_2x_i\otimes y_2)= \psi(a_1x_1\otimes y_1)+ \psi(a_1x_2\otimes y_2)$. Since $a_jU_j=\sup\{a_jx\otimes y: x\otimes y\overset{*}{\leq}U_j\}$,
$\psi(a_iU_i)=f(a_i)W(a_i,U_i)\bot f(a_j)W(a_j,U_j)= \psi(a_jU_j)$ for any $i\not=j$. The last equality follows from the fact that $A=\sup\{a_jU_j:  j\in J\}$.

Let $\{e_{\gamma}: \gamma \in \Gamma \}$ be an orthogonal  basis of $\mathcal H$.  For any finite subset $\{e_j:1\leq j\leq n\}\subset\{e_{\gamma}: \gamma\in\Gamma\}$,  $P_n$  is the projection on  $\bigvee\{e_j:1\leq j\leq n\}$.

  {\bf Claim 3}
  For any $n\geq 1$   and  a finite  rank operators $A$ such that $A\bot I-P_n$, $\varphi(A)\bot I-P_n$.

Take  $a >\|A\|$ and put  $A_{a}=A\oplus a (I-P_n)$.  Let $A=\oplus_{i=1}^k a_i U_i$  be the Penrose   decomposition of $A$. Then $A_a=\sum\limits_{j=1}^ka_jU_j+a( I-P_n)$ is the Penrose decomposition of $A_a$.    Thus $\varphi(A)\bot I-P_n$ by Claim 2.

It follows that $\psi(A)\in P_n\mathcal B(\mathcal H)P_n$ for any $A\in P_n\mathcal B(\mathcal H)P_n$. We now may restrict $\psi$ to $  P_n\mathcal B(\mathcal H)P_n $. This is just the matrix case.  In fact, $f(a)=|g(a)|$ for all $a\in(0,\infty)$ as above.
By \cite[Section 8]{le}(cf. Corollary 8.6),
  the following one  assertion holds.

(1)$\psi(A)=\sum_{i=1}^k g(a_i)U_i$ for any $A=\sum_{i=1}^k a_iU_i\in P_n\mathcal B(\mathcal H)P_n$;

(2)$\psi(A)=\sum_{i=1}^k \overline{g(a_i)}U_i^*$ for any $A=\sum_{i=1}^k a_iU_i\in P_n\mathcal B(\mathcal H)P_n$.

Note that $P_n\mathcal B(\mathcal H)P_n\subseteq P_{n+1}\mathcal B(\mathcal H)P_{n+1}$. If (1) holds for some  $n$, it is elementary  that (1) holds for $n+1$(cf. \cite[Corollary 7.14]{le}).
This means that $(1)$ holds for all finite rank operators $A$.
 Thus (1) holds for all $A\in\mathcal B(\mathcal H)_I$.   Since $\{g(a_j): 1\leq j\leq k\} \subseteq \sigma(|\psi(A)|)$ for any $A=\sum_{i=1}^k a_iU_i\in \mathcal B(\mathcal H)_I$, $g$ is bounded on any bounded subset of $(0,\infty)$. In this case, we define $h(a)=g(a)$, $\forall a\in (0,\infty)$.  Put $S=WU$, $T=U^*$, we have $\varphi(A)=\alpha S( \sum\limits_{j\in J}h(a_j)V_j))T$ for all $A\in\mathcal B(\mathcal H)_I$. Thus $(1)$  of the theorem holds.

 If $(2)$ holds, we define $h(a)=\overline{g(a)}$, $\forall a\in (0,\infty)$.
 In this case, $(2)$ of the theorem holds.
\end{proof}

If  $A=\sum\limits_{j\in J}a_jU_j\in\mathcal B(\mathcal H)_I$ is nonzero, then $|A|=\sum\limits_{j\in J}a_jU_j^*U_j$ and $W_A=\sum\limits_{j\in J}U_j$. For a continuous function $h$ on $(0,\infty)$, we have
 $h(|A|)=\sum\limits_{j\in J}h(a_j)U_j^*U_j$. Moreover,
 $$\sum\limits_{j\in J}h(a_j)U_j=W_Ah(|A|) \mbox{ and }  \sum\limits_{j\in J}h(a_j)U_j^*=W_{A^*}h(|A^*|) .$$
  However, if $h$ is discontinuous, these are fails. In fact, we may construct a bijective function $h$ on $(0,\infty)$ which is not Lebesgue measurable.

\vskip12pt
{\bf Corollary 3.8}
 Let $\varphi$ be a continuous star order automorphism on $\mathcal B(\mathcal H)$. Then there are  a  continuous  complex  function  $h$ on  $[0,\infty)$   satisfying $|h|$ is bijective with $h(0)=0$, a nonzero constant $\alpha\in (0,\infty)$ and two unitary or two anti-unitary  operators  $S$ and $T$ on $\mathcal H$  such that

 $(1)$ $\varphi(A)=\alpha S W_Ah(|A|)T$, $\forall A=W_A|A|\in\mathcal B(\mathcal H)$;

 $(2)$ $\varphi(A)=\alpha S W_{A^*}h(|A^*|)T$, $\forall A=W_A|A|\in\mathcal B(\mathcal H)$.
\vskip12pt

 \begin{proof} We may assume that  $(1)$ holds in Theorem 3.7  for $\varphi|_{\mathcal B(\mathcal H)_I}$.
 If $\varphi$ is continuous, then $h$ is continuous by Theorem 3.7. It is clear $\lim\limits_{x\to 0}h(x)=0$. We may define $h(0)=0$.  For any $A\in\mathcal B(\mathcal H)$, let $A=W_A|A|$ be the polar decomposition of $A$ and $|A|=\int_{\sigma(|A|} x dE_{x}$ the spectral decomposition of $|A|$. Assume that $\sigma (|A|)\subseteq [a,b]$. For any $a=a_0<a_1<a_2<\cdots<a_n=b$,
 Put $E_1=E[a_0,a_1]$ and  $E_j=E(a_{j-1},a_j]$ for $j=2,\cdots n$. Define
 $D_n=\sum\limits_{j=1}^na_jE_j$, $U_j\xi=W_A\xi $  for any $\xi\in E_j\mathcal H$ and $U_j\xi=0$ for any $\xi\in (E_j\mathcal H)^{\bot}$, $1\leq j\leq n$. Then $A_n=W_AD_n=\sum\limits_{j=1}^na_jU_j\in \mathcal B(\mathcal H)_I$ and $A_n\to A(n\to\infty)$. Note that $|A_n|=D_n\to |A|(n\to \infty)$. Thus $\varphi(A)=\lim\limits_{n\to \infty}\varphi(A_n)=\alpha S(\sum\limits_{j=1}^nU_jh(a_j))T=\lim\limits_{n\to\infty}\alpha SW_Ah(D_n)T=\alpha SW_Ah(|A|)T$.

  Similarly we may obtain $(2)$ if $(2)$ holds in Theorem 3.7.
 \end{proof}

 We  now determined all continuous star order automorphisms on $\mathcal B(\mathcal H)$. It is also known  that
 there are discontinuous star order automorphisms on $\mathcal B(\mathcal H)$ by Example 3.2. It may be interesting  to determine all (even continuous) star order automorphisms on $\mathcal B(\mathcal H)_{II}$. On the other hand, the set $PI(\mathcal H)$  of all partial isometries on $\mathcal H$ is a very important poset with respect to star order. We also want to determine star order automorphisms on $PI(\mathcal H)$. If  a star order automorphism on $PI(\mathcal H)$ preserves orthogonality as well, the characterizations  were determined in \cite{sh}.

%\section*{References}

\end{document}